\documentclass{article}

\usepackage{arxiv}

\usepackage[utf8]{inputenc} 
\usepackage[T1]{fontenc}    
\usepackage{hyperref}       
\usepackage{url}            
\usepackage{booktabs}       
\usepackage{amsfonts}       
\usepackage{nicefrac}       
\usepackage{microtype}      
\usepackage{lipsum}		
\usepackage{graphicx}
\usepackage{natbib}
\usepackage{doi}
\usepackage{multirow}
\usepackage{subfigure}
\usepackage{amsmath}

\title{A rotated characteristic decomposition technique for high-order reconstructions in multi-dimensions}

\author{Hua Shen \thanks{This work was supported by the National Natural Science Foundation of China (Contract No. 11901602).} \\
        School of Mathematical Sciences\\
        University of Electronic Science and Technology of China\\
        Chengdu, Sichuan, 611731, China \\
	\texttt{huashen@uestc.edu.cn} \\
	\And
    Matteo Parsani \\
	Extreme Computing Research Center (ECRC)\\
    Computer, Electrical and Mathematical Sciences \& Engineering (CEMSE)\\
    King Abdullah University of Science and Technology (KAUST) \\
    Thuwal, 23955-6900, Kingdom of Saudi Arabia\\
	\texttt{matteo.parsani@kaust.edu.sa} \\
}

\renewcommand{\shorttitle}{\textit{arXiv} Template}

\hypersetup{
pdftitle={A rotated characteristic decomposition technique},
pdfauthor={Hua Shen, Matteo Parsani},
pdfkeywords={characteristic decomposition, high-order schemes, hyperbolic conservation laws, WENO, finite volume},
}

\begin{document}
\maketitle

\begin{abstract}
When constructing high-order schemes for solving hyperbolic conservation laws,
the corresponding high-order reconstructions are commonly performed in characteristic spaces 
to eliminate spurious oscillations as much as possible.
For multi-dimensional finite volume (FV) schemes,
we need to perform the characteristic decomposition several times in different normal directions of the target cell, 
which is very time-consuming.
In this paper, we propose a rotated characteristic decomposition technique 
which requires only one-time decomposition for multi-dimensional reconstructions. 
The rotated direction depends only on the gradient of a specific physical quantity which is cheap to calculate.
This technique not only reduces the computational cost remarkably, but also controls spurious oscillations effectively.
We take a third-order weighted essentially non-oscillatory finite volume (WENO-FV) scheme for solving the Euler equations
as an example to demonstrate the efficiency of the proposed technique.
\end{abstract}


\section{Introduction}
\label{SEC:Intro}
The solution to a nonlinear hyperbolic conservation law may become discontinuous even when the initial conditions are sufficiently smooth.
Therefore, the shock-capturing property is a desired feature of a numerical scheme for solving hyperbolic conservation laws.
Usually, low-order schemes are more stable than high-order schemes for capturing discontinuities.
However, in smooth regions, low-order schemes converge too slow, so they require a very fine mesh to achieve low-level errors.
In many applications, discontinuities and sophisticated structures coexist, hence shock-capturing high-order schemes are preferred.
In the development history of numerical schemes for solving hyperbolic conservation laws,
the order of the numerical schemes gradually increased.
Van Leer \cite{vanLeer1979MUSCLV} proposed a second-order monotonic upstream-centered scheme for conservation laws
(MUSCL) which is an extension of the first-order Godunov scheme \cite{Godunov1959}.
Colella and Woodward \cite{Colella1984PPM} proposed a third-order piecewise parabolic method (PPM) which was extended from MUSCL.
Harten and his coworkers \cite{Harten1987ENO} proposed a general framework to design uniformly high-order accurate essentially non-oscillatory (ENO) schemes by using the strategy of adaptively choosing the smoothest stencil from several local candidates.
Liu \emph{et al}. \cite{Liu1994WENO} proposed weighted ENO (WENO) schemes by assigning a proper weight to each candidate stencil.
Jiang and Shu \cite{Jiang_Shu1996WENO} improved the accuracy of WENO schemes by carefully design the weights.
Borges \emph{et al}. \cite{Borges2008WENO_Z} proposed WENO-Z schemes by further improving the weights for the candidate stencils.
Levy, Puppo and Russo \cite{Levy1999CWENO, Levy2000CWENO_ANM, Levy2000CompactCWENO} 
proposed a class of central WENO (CWENO) schemes 
which reconstruct a single polynomial in the entire traget cell.
Therefore, CWENO schemes are cheaper than standard WENO-FV schemes
when there are many reconstruction points in a cell 
and can be trivially extended to the case of non-Cartesian meshes. 
A detailed literature review of CWENO schemes can be found in \cite{Cravero2019CWENO}.
In principle, we can construct arbitrarily high-order ENO/WENO schemes, but the size of stencils increases as the order increases.
This feature makes it unfriendly for their implementation on unstructured meshes, although we can do it,
see for example \cite{Abgrall1994ENO_Tri, Hu1999WENO_Tri, Dumbser2007WENO_FV}.
To overcome the above shortcoming, some compact high-order schemes are constructed by increasing the internal degree of freedoms (DOFs).
The representative of compact schemes is the discontinuous Galerkin (DG) scheme \cite{Reed1973DG, Cockburn1989RKDGII, Cockburn1989RKDGIII, Cockburn1990RKDGIV, Cockburn1998RKDGV} that relies on nonlinear limiters to capture discontinuities,
see for example \cite{Qiu2004HermiteWENO, Qiu2005HermiteWENO2D, Zhong2013DG_Limiter, Zhu2016DG_Limiter}.
Dumbser \emph{et al}. \cite{Dumbser2008PNPM, Dumbser2009PNPM, Dumbser2010PNPM}
proposed a class of $P_NP_M$ schemes which can flexibly choose the internal DOFs and reconstruction DOFs,
thereby adjusting the compactness of the stencil.

When solving a hyperbolic system by ENO/WENO schemes or DG schemes with non-linear limiters,
it is necessary to carry out reconstructions in characteristic space to control spurious oscillations \cite{Qiu2002Decomposition}.
The reconstructions of multi-dimensional ENO/WENO finite difference (ENO/WENO-FD) schemes \cite{Shu1988EfficientENO, Jiang_Shu1996WENO}
can be efficiently performed in a dimension-by-dimension manner.
However, ENO/WENO-FD only works for regular meshes.
For irregular meshes, we have to resort to FV-type schemes,
but classical ENO/WENO-FV schemes and DG schemes need to
perform characteristic decomposition and reconstructions in all normal directions of the target cell.
This is one of the major reasons that makes ENO/WENO-FV schemes more expensive than ENO/WENO-FD schemes in multi-dimensions.
For Cartesian meshes, there are two and three normal directions in two- and three-dimensions, respectively.
For unstructured meshes, there are even more normal directions than Cartesian meshes.
In order to reduce the computational cost, 
some scholars proposed to only perform characteristic decomposition near discontinuities 
\cite{Ren2003hybridCompatWENO, Puppo2003adaptive, Puppo2011adaptive, Li2010hybridWENO, Peng2019adaptiveWENOZ}.
This approach relies on a good switch function, 
otherwise the transition between non-characteristic reconstructions and characteristic reconstructions may become nonsmooth.
In this paper, we introduce a rotated characteristic decomposition technique
which requires only one-time characteristic decomposition in the direction of a physical quantity's gradient.
This technique not only reduces the amount of computation significantly, but also can eliminate spurious oscillations effectively.

\section{Description of the rotated characteristic decomposition technique}
\label{SEC:RCD}
\subsection{A brief review of WENO-FV schemes}
\label{SubSEC:WENO_FV}
We consider the hyperbolic conservation law, equipped with certain initial conditions and boundary conditions,
which can be expressed as
\begin{subequations}\label{Eq:HCL}
  \begin{equation}\label{SubEq:HCL_Eq}
    \frac{\partial \mathbf{U}}{\partial t}+\nabla\cdot\mathbf{F}=0,\quad \mathbf{x}\in\Omega,\quad t\in[0,\infty),
  \end{equation}
  \begin{equation}\label{SubEq:HCL_IC}
    \mathbf{U}(\mathbf{x},0)=\mathbf{U}_0(\mathbf{x}), \quad \mathbf{x}\in\Omega,
  \end{equation}
  \begin{equation}\label{SubEq:HCL_BC}
    \mathbf{U}(\mathbf{x},t)=\mathbf{U}_{\partial\Omega}(\mathbf{x},t), \quad \mathbf{x}\in\partial\Omega, \quad t\in[0,\infty).
  \end{equation}
\end{subequations}
A semi-discrete FV scheme for solving Eq. (\ref{Eq:HCL}) can be written as
\begin{equation}\label{Eq:SemiDisFV}
  \frac{d\bar{\mathbf{U}_i}}{dt}=-\frac{1}{\Omega_i}\sum_{j=1}^{J}|\partial \Omega_{i,j}|
  \sum_{k=1}^{K}\omega_k\mathbf{\hat{F}}\left(\mathbf{x}_k^G\right)\cdot\mathbf{n}_{i,j},
\end{equation}
where $\Omega_i$ is the discrete cell, 
$\partial\Omega_{i,j}$ is the $j$th edge of $\Omega_i$,
$\mathbf{n}_{i,j}$ is the unit outward normal vector of $\partial \Omega_{i,j}$,
$\mathbf{x}_k^G$ and $\omega_k$  are the Gaussian quadrature points and weights
applied for the computation of the numerical flux on $\partial\Omega_{i,j}$.
The numerical fluxes at the Gaussian quadrature points can be calculated by local Riemann solvers.
Here, we use the local Lax–Friedrichs flux that is expressed as
\begin{equation}\label{Eq:LxFFlux}
  \mathbf{\hat{F}}\left(\mathbf{x}_k^G\right)=\frac{1}{2}\left\{\mathbf{F}^+\left(\mathbf{x}_k^G\right)+\mathbf{F}^-\left(\mathbf{x}_k^G\right)
  -S_{i,j}\left[\mathbf{U}^+\left(\mathbf{x}_k^G\right)-\mathbf{U}^-\left(\mathbf{x}_k^G\right)\right]\right\},
\end{equation}
where $S_{i,j}$ is the local upper bound for the eigenvalues of the Jacobian in the $\mathbf{n}_{i,j}$ direction,
$\mathbf{U}^\pm\left(\mathbf{x}_k^G\right)$ is interpolated from the target cell and the neighboring cell.

In order to get a $p$th-order scheme, we have to reconstruct a $(p-1)$th-order polynomial $\mathbf{U}_i(\mathbf{x})$ from cell averages
to approximate the solution in every cell.
A popular reconstruction method is WENO which is a nonlinearly convex combination of linear reconstructions on several local stencils.
In this work, we take a two-dimensional third-order WENO-FV scheme on Cartesian meshes
as an example to demonstrate the efficiency of the proposed technique.
The WENO reconstructions proposed by Balsara \emph{et al}. \cite{Balsara2009ADER-WENO} is used for space reconstructions,
and the third order TVD Runge–Kutta method \cite{Shu1988EfficientENO} is adopted for time discretization.
When constructing non-linear stencil weights, we use the smoothness measures to the second power and set $\epsilon=10^{-40}$.
The other settings exactly keep the same with those used by Balsara \emph{et al}. \cite{Balsara2009ADER-WENO}.
\subsection{The rotated characteristic decomposition technique for the two-dimensional Euler equations}
\label{SubSEC:Euler_RCD}
We consider the two-dimensional Euler equations, namely
\begin{subequations}\label{Eq:2D_Euler}
  \begin{equation}\label{SubEq:2D_Euler_1}
    \frac{\partial \mathbf{U}}{\partial t}+\frac{\partial \mathbf{F}_1}{\partial x_1}+\frac{\partial \mathbf{F}_2}{\partial x_2}=0,
  \end{equation}
  \text {with}
  \begin{equation}\label{SubEq:2D_Euler_2}
    \mathbf{U}=\begin{bmatrix}
                 \rho \\
                 \rho u \\
                 \rho v\\
                 \rho e
               \end{bmatrix},
    \mathbf{F}_1=\begin{bmatrix}
                 \rho u \\
                 \rho u^2+p \\
                 \rho uv\\
                 (\rho e+p)u
               \end{bmatrix},
    \mathbf{F}_2=\begin{bmatrix}
                 \rho v \\
                 \rho uv \\
                 \rho v^2+p\\
                 (\rho e+p)v
               \end{bmatrix},
  \end{equation}
\end{subequations}
where $\rho$ is the density, $u$ and $v$ are respectively velocity component in $x_1$ and $x_2$ directions, $p$ is the pressure,
$e$ is the specific total energy.
For ideal gases, we have
\begin{equation}\label{Eq:TotalE}
  e=\frac{p}{(\gamma-1)\rho}+\frac{1}{2}(u^2+v^2).
\end{equation}

The Jacobian matrix corresponding to $\mathbf{F}_1$ is
\begin{equation}\label{Eq:DF1DU}
  \mathbf{A}(\mathbf{U})=\frac{\partial \mathbf{F}_1}{\partial \mathbf{U}}=
  \begin{bmatrix}
    0                                   & 1               & 0            & 0 \\
    \frac{\gamma-1}{2}K-u^2             & (3-\gamma)u     & (1-\gamma)v  & \gamma-1 \\
    -uv                                 & v               & u            & 0 \\
    \left(\frac{\gamma-1}{2}K-H\right)u & H-(\gamma-1)u^2 & (1-\gamma)uv & \gamma u
  \end{bmatrix},
\end{equation}
where $K=u^2+v^2$, and $H=e+\frac{p}{\rho}$.
The eigenvalues of $\mathbf{A}(\mathbf{U})$ are
\begin{subequations}\label{Eq:2DEulerEigen}
  \begin{equation}\label{Eq:2DEulerEigenValue}
    \lambda_1=u-a,\quad \lambda_2=\lambda_3=u,\quad \lambda_3=u+a,
  \end{equation}
\text{of which the corresponding right eigenvectors are}
  \begin{equation}\label{Eq:2DEulerEigenVector}
    \mathbf{r}_1=\begin{bmatrix}
                   1 \\
                   u-a \\
                   v \\
                   H-au
                 \end{bmatrix},
    \mathbf{r}_2=\begin{bmatrix}
                   1 \\
                   u \\
                   v \\
                   \frac{K}{2}
                 \end{bmatrix},
    \mathbf{r}_3=\begin{bmatrix}
                   0 \\
                   0 \\
                   1 \\
                   v
                 \end{bmatrix},
    \mathbf{r}_4=\begin{bmatrix}
                   1 \\
                   u+a \\
                   v \\
                   H+au
                 \end{bmatrix},
  \end{equation}
\end{subequations}
where $a=\sqrt{\frac{\gamma p}{\rho}}$ is the sound speed.
Define the right eigenvector matrix as $\mathbf{R}=[\mathbf{r}_1, \mathbf{r}_2, \mathbf{r}_3, \mathbf{r}_4]$,
then the Jacobian matrix can be written as $\mathbf{A}(\mathbf{U})=\mathbf{R}\mathbf{\Lambda}\mathbf{R}^{-1}$,
where $\mathbf{\Lambda}$ is the diagonal matrix with $\mathbf{\Lambda}_{kk}=\lambda_k$,
the left eigenvector matrix $\mathbf{R}^{-1}$ is expressed as
\begin{equation}\label{Eq:2DEuler_LeftEigenV}
  \mathbf{R}^{-1}
  =\frac{\gamma-1}{2a^2}\begin{bmatrix}
  \frac{au}{\gamma-1}+\frac{K}{2} & -\frac{a}{\gamma-1}-u & -v & 1 \\
  \frac{2a^2}{\gamma-1}-K & 2u & 2v & -2 \\
  -\frac{2a^2v}{\gamma-1} & 0 & \frac{2a^2}{\gamma-1} & 0 \\
  -\frac{au}{\gamma-1}+\frac{K}{2} & \frac{a}{\gamma-1}-u & -v & 1
                  \end{bmatrix}.
\end{equation}

Once we get the flux and its eigenstructure in $x_1$ direction, we can use the rotational invariance property of the Euler equations
to calculate the flux and the corresponding eigenstructure in an arbitrary direction, $\mathbf{n}=(n_1,n_2)=(cos\theta,sin\theta)$, as \cite{Toro2013CFDBook}
\begin{subequations}\label{Eq:RotatedEigen}
  \begin{equation}\label{Eq:RotatedFlux}
    \mathbf{F}(\mathbf{U},\theta)=cos\theta\mathbf{F}_1+sin\theta\mathbf{F}_2=\mathbf{T}^{-1}\mathbf{F}_1\left(\mathbf{T}\mathbf{U}\right),
  \end{equation}
  \begin{equation}\label{Eq:RotatedJacobian}
    \mathbf{A}(\mathbf{U},\theta)=cos\theta\frac{\partial \mathbf{F}_1}{\partial \mathbf{U}}+sin\theta\frac{\partial \mathbf{F}_2}{\partial \mathbf{U}}
    =\mathbf{T}^{-1}\mathbf{A}\left(\mathbf{T}\mathbf{U}\right)\mathbf{T},
  \end{equation}
  \begin{equation}\label{Eq:RotatedEigenV}
    \mathbf{R}(\mathbf{U},\theta)=\mathbf{T}^{-1}\mathbf{R}\left(\mathbf{T}\mathbf{U}\right),\quad
     \mathbf{R}^{-1}(\mathbf{U},\theta)=\mathbf{R}^{-1}\left(\mathbf{T}\mathbf{U}\right)\mathbf{T},
  \end{equation}
\text{with the rotation matrix and its inverse defined as}
  \begin{equation}\label{Eq:RotatedMatrix}
    \mathbf{T}=\begin{bmatrix}
                 1 & 0 & 0 & 0 \\
                 0 & cos\theta & sin\theta & 0 \\
                 0 & -sin\theta & cos\theta & 0 \\
                 0 & 0 & 0 & 1
               \end{bmatrix},
    \mathbf{T}^{-1}=\begin{bmatrix}
                 1 & 0 & 0 & 0 \\
                 0 & cos\theta & -sin\theta & 0 \\
                 0 & sin\theta & cos\theta & 0 \\
                 0 & 0 & 0 & 1
               \end{bmatrix}.
  \end{equation}
\end{subequations}

When we apply WENO schemes to solve the Euler equations,
we have to perform all the reconstructions in characteristic space to eliminate spurious oscillations as much as possible.
In standard WENO approaches, before we calculate the flux in $\mathbf{n}_{i,j}=(cos\theta_{i,j},sin\theta_{i,j})$ direction in Eq. (\ref{Eq:SemiDisFV}),
we need the following three major steps in advance:\\
{\bfseries Step 1}. Project the cell averages $\bar{\mathbf{U}}_m$ within the reconstruction stencils of the target cell $\Omega_i$ to local characteristic space, i.e., $\bar{\mathbf{U}}_{m,j}^{(c)}=\mathbf{R}^{-1}(\bar{\mathbf{U}}_i,\theta_{i,j})\bar{\mathbf{U}}_m$;\\
{\bfseries Step 2}. Perform WENO reconstruction based on $\bar{\mathbf{U}}_{m,j}^{(c)}$,
i.e., $\mathbf{U}_{i,j}^{(c)}(\mathbf{x})=WENOREC\left(\bar{\mathbf{U}}_{m,j}^{(c)}\right)$;\\
{\bfseries Step 3}. Transform $\mathbf{U}_{i,j}^{(c)}(\mathbf{x})$ back to physical space, i.e.,
$\mathbf{U}_{i,j}(\mathbf{x})=\mathbf{R}(\bar{\mathbf{U}}_i,\theta_{i,j})\mathbf{U}_{i,j}^{(c)}(\mathbf{x})$.\\
Note that, we need to repeat $N$ times the above steps on a target cell with $N$ different normal directions.
For example, on two-dimensional Cartesian meshes, we need to execute the above steps in $x_1$ and $x_2$ respectively.

In order to improve the efficiency, we propose a rotated characteristic decomposition technique which contains the following steps:\\
{\bfseries Step 1}. Determine the rotated direction as $\hat{\mathbf{n}}_i=(cos\hat{\theta}_i,sin\hat{\theta}_i)$ with
\begin{equation}\label{Eq:Gradient}
  cos\hat{\theta}_i=\frac{\rho_{x_1,i}+\epsilon}{\sqrt{(\rho_{x_1,i}+\epsilon)^2+(\rho_{x_2,i}+\epsilon)^2}}, 
  sin\hat{\theta}_i=\frac{\rho_{x_2,i}+\epsilon}{\sqrt{(\rho_{x_1,i}+\epsilon)^2+(\rho_{x_2,i}+\epsilon)^2}},
\end{equation}
where the spatial derivatives $\rho_{x_1,i}$ and $\rho_{x_2,i}$ are calculated by central differences,
and $\epsilon=10^{-40}$ is a small number used to avoid singularity;\\
{\bfseries Step 2}. Project the cell averages $\bar{\mathbf{U}}_m$ within the reconstruction stencils of the target cell $\Omega_i$ to local characteristic space, i.e., $\bar{\mathbf{U}}_m^{(c)}=\mathbf{R}^{-1}(\bar{\mathbf{U}}_i,\hat{\theta}_i)\bar{\mathbf{U}}_m$;\\
{\bfseries Step 3}. Perform WENO reconstruction based on $\bar{\mathbf{U}}_m^{(c)}$,
i.e., $\mathbf{U}_i^{(c)}(\mathbf{x})=WENOREC\left(\bar{\mathbf{U}}_m^{(c)}\right)$;\\
{\bfseries Step 4}. Transform $\mathbf{U}_i^{(c)}(\mathbf{x})$ back to physical space, i.e.,
$\mathbf{U}_i(\mathbf{x})=\mathbf{R}(\bar{\mathbf{U}}_i,\hat{\theta}_i)\mathbf{U}_i^{(c)}(\mathbf{x})$.\\
At the first glance, the rotated characteristic decomposition technique requires one more step than the standard technique.
However, the additional step to determine the rotated direction is very cheap,
and the rotated characteristic decomposition technique always requires only one-time characteristic decomposition
regardless of the number of the normal directions of the target cell.
Therefore, the rotated characteristic decomposition technique can significantly reduce the computational cost comparing to the standard technique.
Furthermore, spurious oscillations usually appear in the direction of shocks and contact discontinuities
which can be detected by the density gradient,
so we can effectively control the spurious oscillations by performing characteristic decomposition in the direction of the density gradient.


%
\section{Numerical examples}
\label{SEC:NumExam}
We serially run all the simulations on a laptop equipped with Intel(R) Core(TM)i7-4720HQ CPU@2.60GHz.
In all simulations, we set $CFL=0.8$ and $\gamma=1.4$.
In this section, we respectively denote no characteristic decomposition as NCD, standard characteristic decomposition as SCD,
and rotated characteristic decomposition as RCD for short.
We note that, we perform RCD in the density gradient direction and the total energy gradient direction respectively,
but no difference was observed.
Therefore, we show only the results of RCD in the density gradient direction.
\subsection{Isentropic vortex evolution}
\label{SubSEC:Vortex}
This case is used to demonstrate the effect of the RCD technique on smooth flows.
The computational domain is $[-5, 5]\times[-5, 5]$  of which all boundaries are periodic.
Initially, the mean flow ($\rho=1, p=1$ and $(u,v)=(1,1)$) is disturbed by an isentropic vortex which is described by
\begin{equation}\label{Eq:Vortex_Perturbations}
  (\delta u, \delta v)=\frac{\psi}{2\pi}e^{0.5(1-r^2)}(-y, x),\quad \delta T=-\frac{(\gamma-1)\psi^2}{8\gamma\pi^2}e^{(1-r^2)},
   \quad \delta S=0,
\end{equation}
where $r^2=x^2+y^2$, $S=p/\rho^\gamma$ is the entropy, $T=p/\rho$ is the temperature, and $\psi$ is the vortex strength.
In our computations, we set $\psi=5$.

The exact solution of this problem can be portrayed by the vortex moving at the mean velocity.
We use the third-order WENO scheme with different characteristic decomposition techniques to compute this problem to $t=2$.
Table \ref{tab:vortex} shows the corresponding density errors, convergence rates, and computational costs.
As we can see, the RCD technique has no effect on the accuracy of the scheme,
but it significantly reduces CPU times as comparing with the SCD technique.

\begin{table}
\caption{The density errors, convergence rates, and computational costs of the third-order WENO-FV scheme
for the isentropic vortex evolution problem.}
\label{tab:vortex}       
\begin{tabular}{cccccc}
\hline
\multicolumn{6}{c}{WENO-FV with NCD}\\
\hline
cell number & $L_1$ & Order & $L_\infty$ &order &CPU time (s) \\
\hline
$50\times50$ & 4.94E-4 &-& 7.74E-3&-& 0.8 \\
$100\times100$ & 6.89E-5 & 2.84&9.32E-4&3.05&5.8 \\
$150\times150$ &2.13E-5  & 2.90&2.34E-4&3.41&19\\
$200\times200$ &9.51E-6  &2.80 &1.00E-4&2.96&45\\
\hline
\multicolumn{6}{c}{WENO-FV with SCD}\\
\hline
cell number & $L_1$ & Order & $L_\infty$ &order &CPU time (s) \\
\hline
$50\times50$ & 4.99E-4 &-& 7.78E-3&-& 2 \\
$100\times100$ & 6.87E-5 & 2.86&9.14E-4&3.09&16 \\
$150\times150$ &2.12E-5  & 2.90&2.29E-4&3.41&52\\
$200\times200$ &9.49E-6  &2.79 &1.01E-4&2.85&123\\
\hline
\multicolumn{6}{c}{WENO-FV with RCD}\\
\hline
cell number & $L_1$ & Order & $L_\infty$ &order &CPU time (s) \\
\hline
$50\times50$ & 5.01E-4 &-& 7.72E-3&-& 1.2 \\
$100\times100$ & 6.87E-5 & 2.87&9.17E-4&3.07&9 \\
$150\times150$ &2.12E-5  & 2.90&2.31E-4&3.40&30\\
$200\times200$ &9.49E-6  &2.79 &1.01E-4&2.88&70\\
\hline
\end{tabular}
\end{table}

\subsection{Double Mach reflection problem}
\label{SubSEC:DMR}
This problem was proposed by Woodward and Colella \cite{Woodward1984JCP}.
The computational domain is $[0, 4]\times[0, 1]$.
At $t=0$, an oblique Mach 10 shock is inclined at an angle of $60^\circ$ to the horizontal direction.
The pre-shock and post-shock states are given by
\begin{equation}\label{Eq:DMR}
\left(\rho, u, v, p\right)=
  \begin{cases}
    \left(1.4, 0, 0, 1\right) & \text{if } x>\frac{1}{6}+\frac{y}{\sqrt{3}}, \\
    \left(8, 8.25sin(60^\circ), -8.25cos(60^\circ), 116.5\right), & \text{otherwise}.
  \end{cases}
\end{equation}
The post-shock states are imposed on the left boundary,
and nonreflective boundary conditions are applied to the right boundary.
On the top, the time-dependent boundary conditions are determined by the exact motion of the oblique shock.
On the bottom, outflow conditions are imposed for $x\le\frac{1}{6}$,
and reflective boundary conditions are imposed for $x>\frac{1}{6}$.

\begin{figure}
  \centering
  \subfigure[WENO-FV with NCD]{
  \label{FIG:DMR_NCD}
  \includegraphics[width=10 cm]{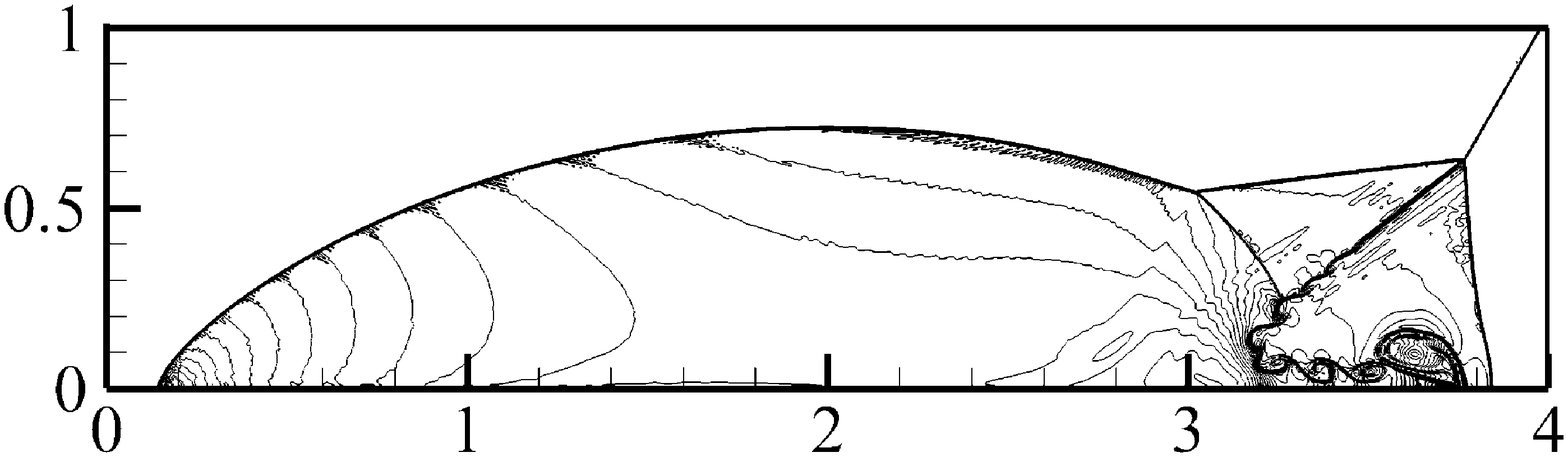}}
  \subfigure[WENO-FV with SCD]{
  \label{FIG:DMR_SCD}
  \includegraphics[width=10 cm]{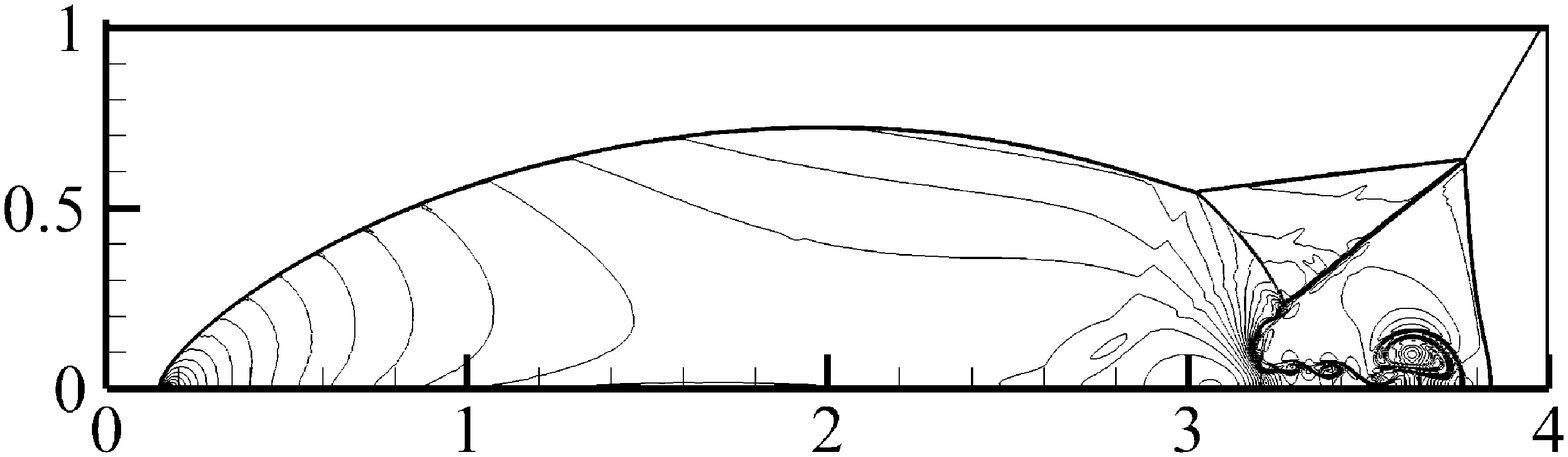}}
  \subfigure[WENO-FV with RCD]{
  \label{FIG:DMR_RCD}
  \includegraphics[width=10 cm]{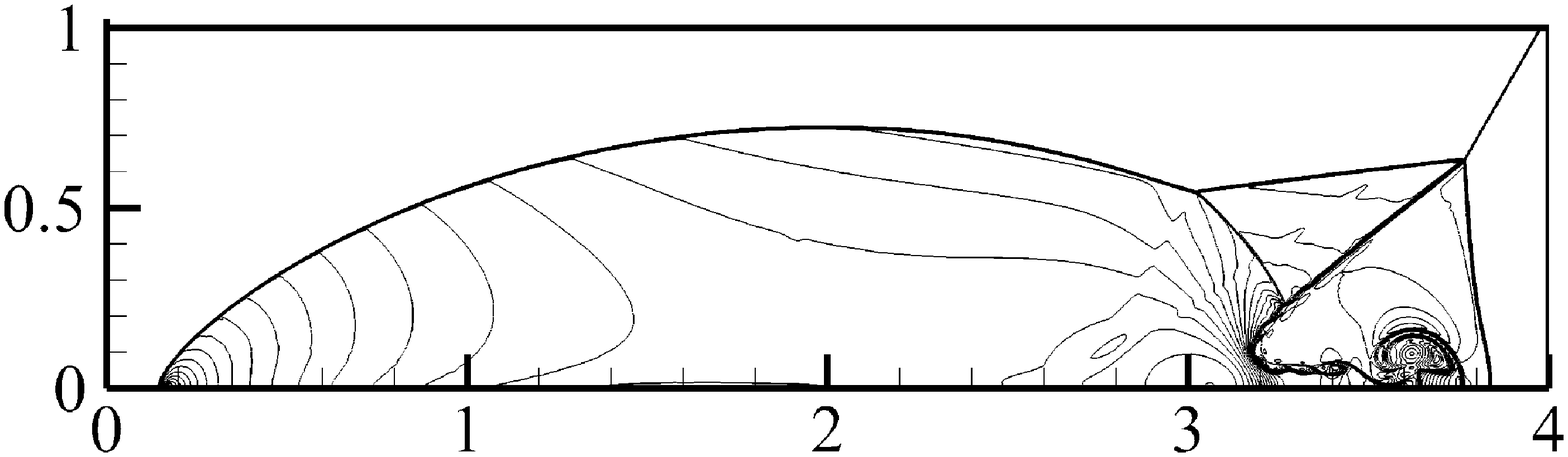}}
  \caption{The entire view of the density contours of the double Mach reflection problem at $t=0.28$
  calculated by the third-order WENO-FV scheme with different characteristic decomposition techniques.
  The mesh size is 1/480 and the density contours contain 50 equidistant contours from 2 to 22.}
\label{FIG:DMR_Density_Contour}
\end{figure}

\begin{figure}
  \centering
  \subfigure[WENO-FV with NCD]{
  \label{FIG:DMR_NCD_enlarge}
  \includegraphics[width=5.5 cm]{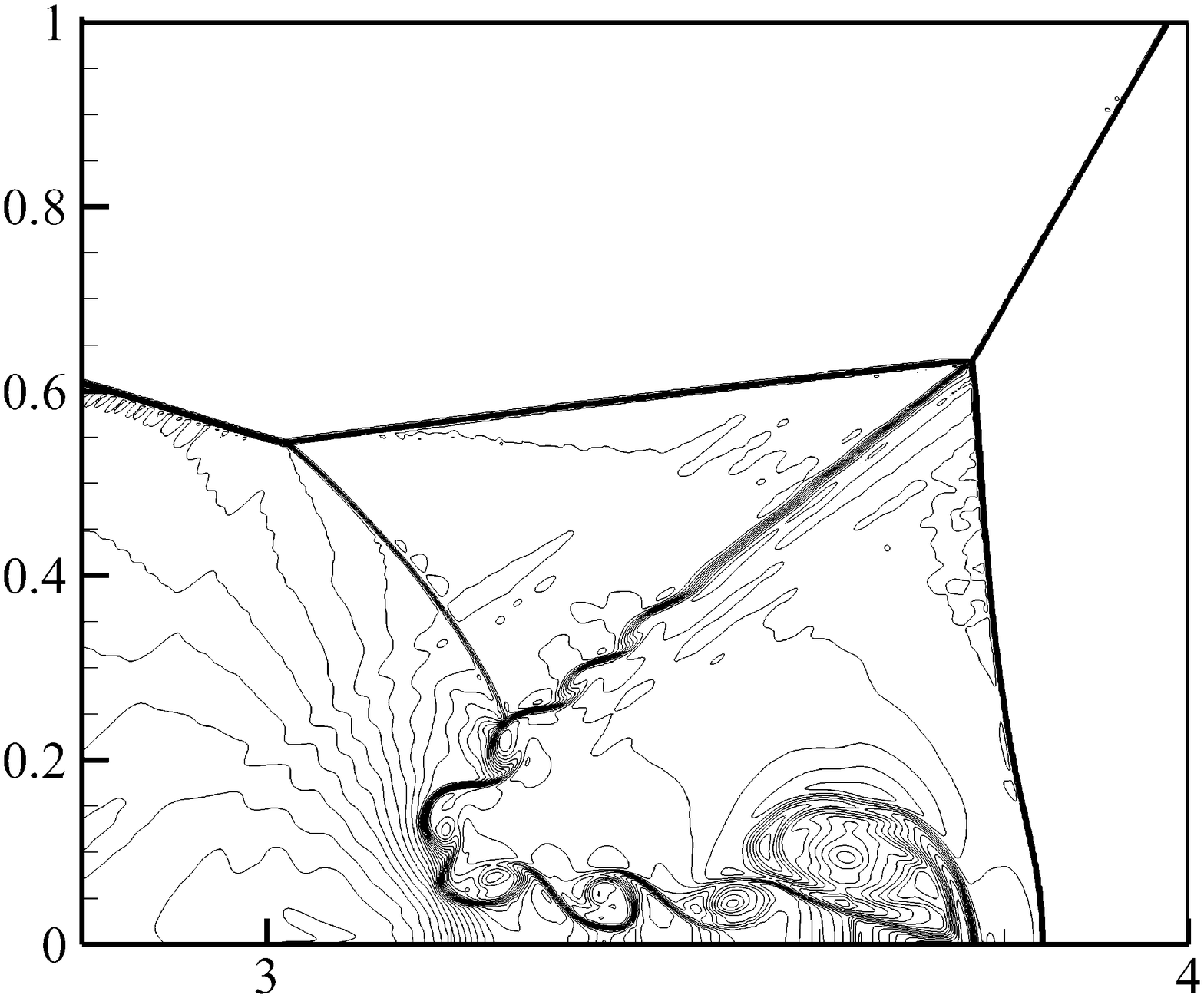}}
  \subfigure[WENO-FV with SCD]{
  \label{FIG:DMR_SCD_enlarge}
  \includegraphics[width=5.5 cm]{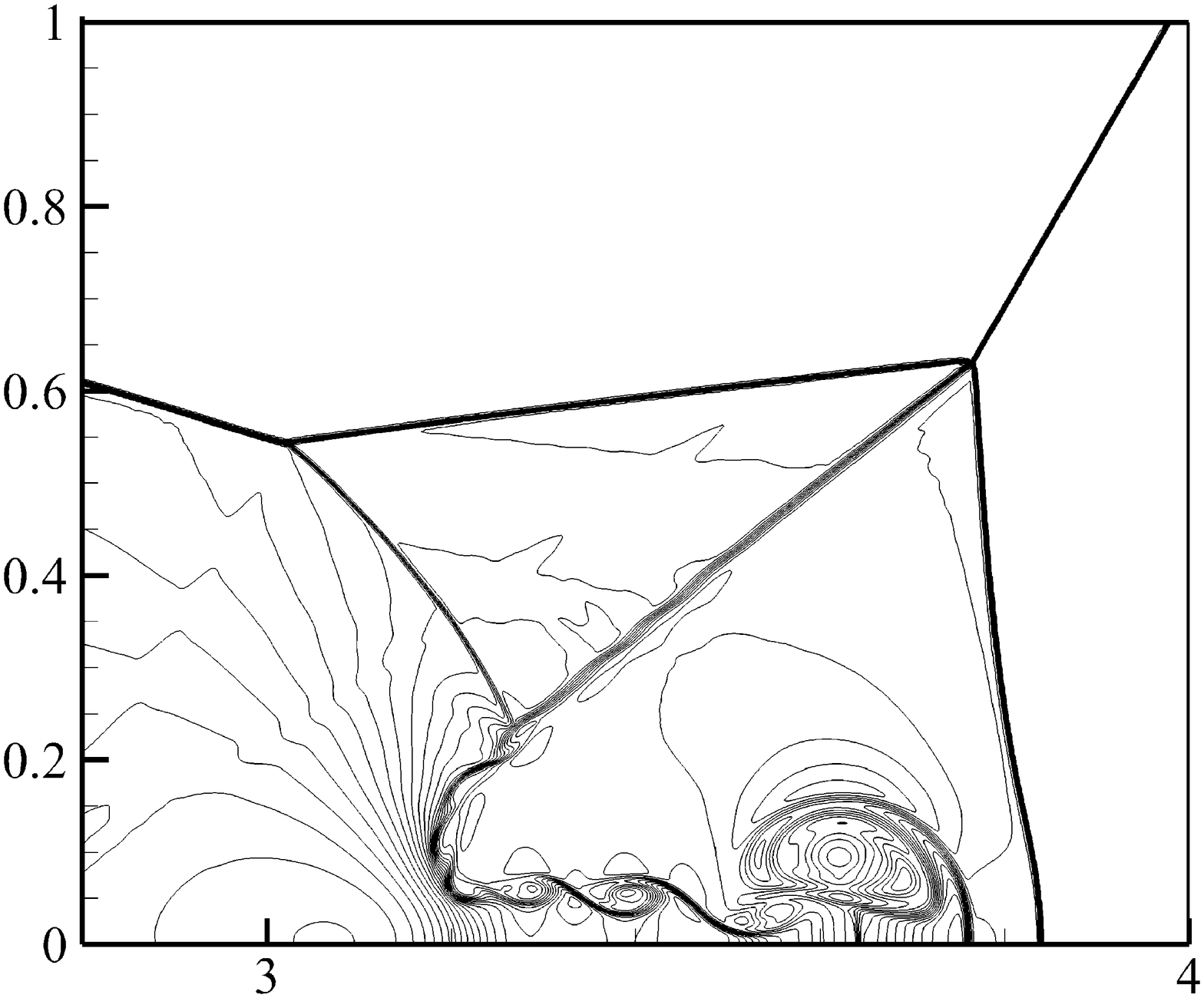}}
  \subfigure[WENO-FV with RCD]{
  \label{FIG:DMR_RCD_enlarge}
  \includegraphics[width=5.5 cm]{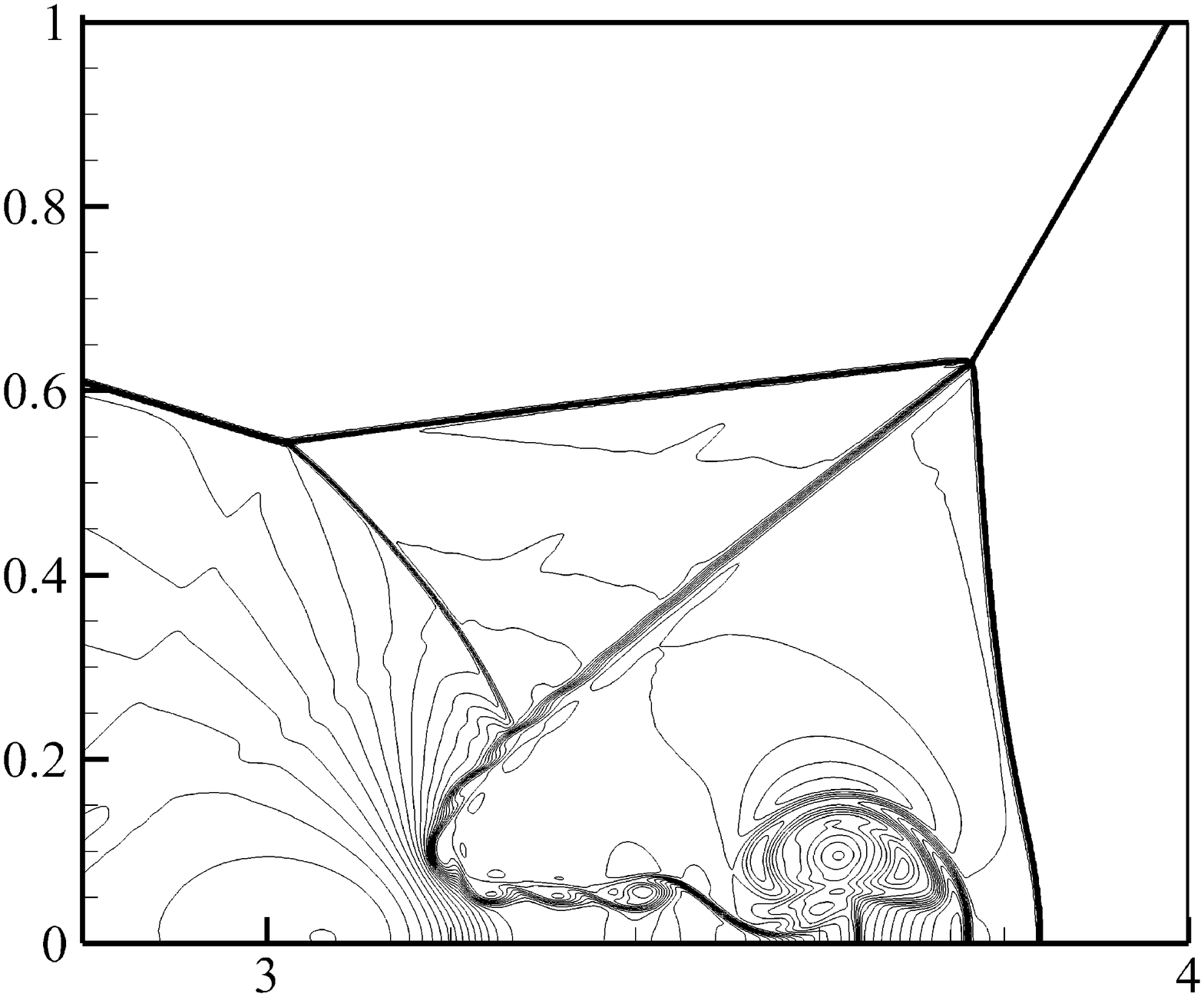}}
  \caption{The enlarged view of the density contours of the double Mach reflection problem at $t=0.28$
  calculated by the third-order WENO-FV scheme with different characteristic decomposition techniques.
  The mesh size is 1/480 and the density contours contain 50 equidistant contours from 2 to 22.}
\label{FIG:DMR_Density_Contour_enlarge}
\end{figure}

This problem is equivalent to a horizontally moving Mach 10 shock reflects by a $60^\circ$ wedge.
Under current condition, a double Mach reflection will occur.
We run this problem by the third-order WENO-FV scheme with $1920\times480$ cells up to $t=0.28$.
The CPU times of WENO-FV with NCD, SCD, and RCD are 18925s, 54032s, and 32179s, respectively.
Fig. \ref{FIG:DMR_Density_Contour} shows the entire view of the density contours,
and Fig. \ref{FIG:DMR_Density_Contour_enlarge} shows the corresponding enlarged view.
As we expected, the WENO-FV with NCD cannot control spurious oscillations effectively.
The RCD technique can control oscillations well with a significantly less CPU time than the SCD technique.
Meanwhile, the RCD and NCD techniques achieve similar resolutions for the complex slip line instability,
as shown by Fig. \ref{FIG:DMR_Density_Contour_enlarge}.
\subsection{Shock reflection on a plate}
\label{SubSEC:shock_ref}
This is a steady-state shock reflection problem which was proposed by Yee \emph{et al}. \cite{Yee1985JCP}.
The computational domain is $[0, 4]\times[0, 1]$ which is initialized by $\left(\rho, u, v, p\right)=\left(1, 2.9, 0, 1/\gamma\right)$.
Non-reflective boundary conditions are imposed on the right, and reflective boundary conditions are imposed on the bottom.
On the left and top, we impose the following Dirichlet boundary conditions
\begin{equation}\label{Eq:DMR}
\left(\rho, u, v, p\right)=
  \begin{cases}
    \left(1, 2.9, 0, 1/\gamma\right) & \text{on the left}, \\
    \left(1.69997, 2.61934, -0.50632, 1.52819\right), & \text{on the top}.
  \end{cases}
\end{equation}

\begin{figure}
  \centering
  \subfigure[WENO-FV with NCD]{
  \label{FIG:Shock_Ref_NCD}
  \includegraphics[width=10 cm]{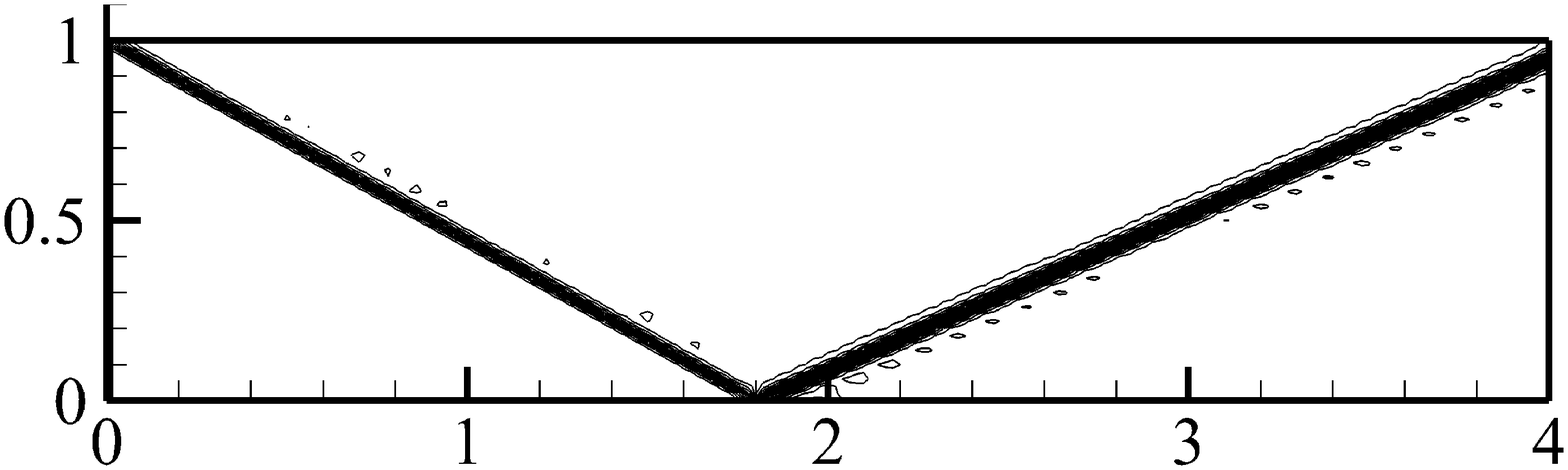}}
  \subfigure[WENO-FV with SCD]{
  \label{FIG:Shock_Ref_SCD}
  \includegraphics[width=10 cm]{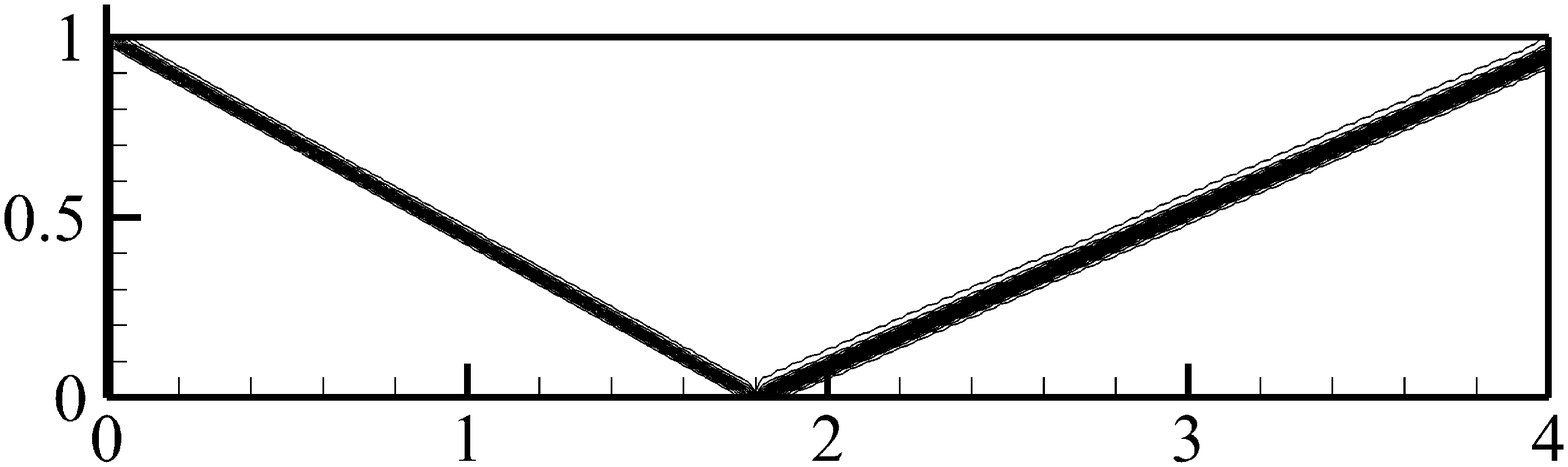}}
  \subfigure[WENO-FV with RCD]{
  \label{FIG:Shock_Ref_RCD}
  \includegraphics[width=10 cm]{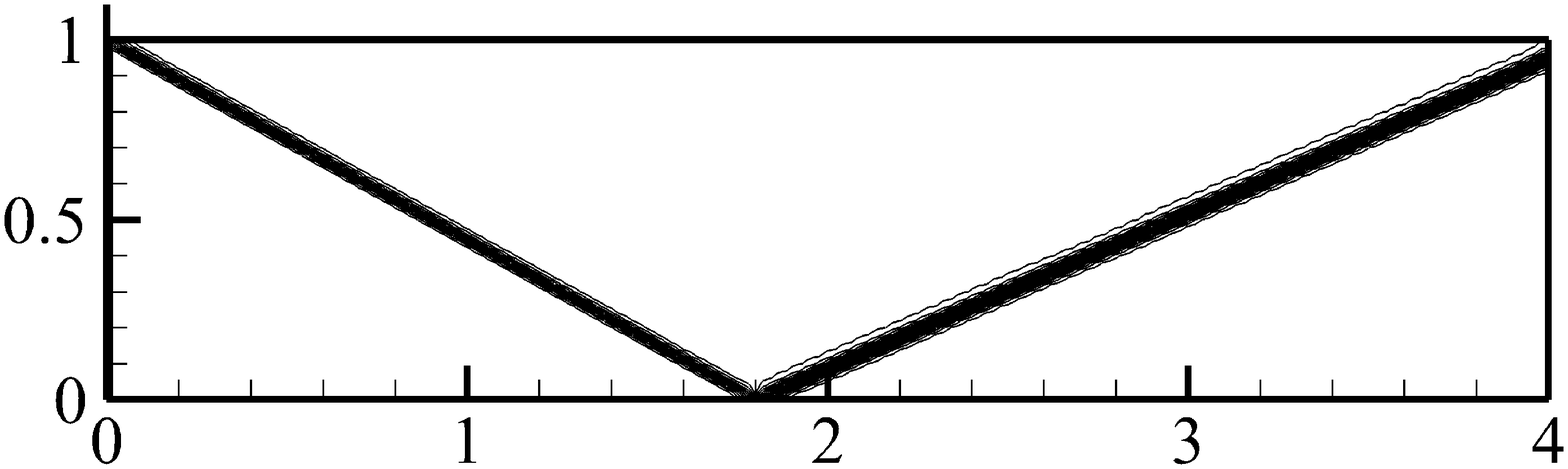}}
  \subfigure[WENO-FV with xCD]{
  \label{FIG:Shock_Ref_xCD}
  \includegraphics[width=10 cm]{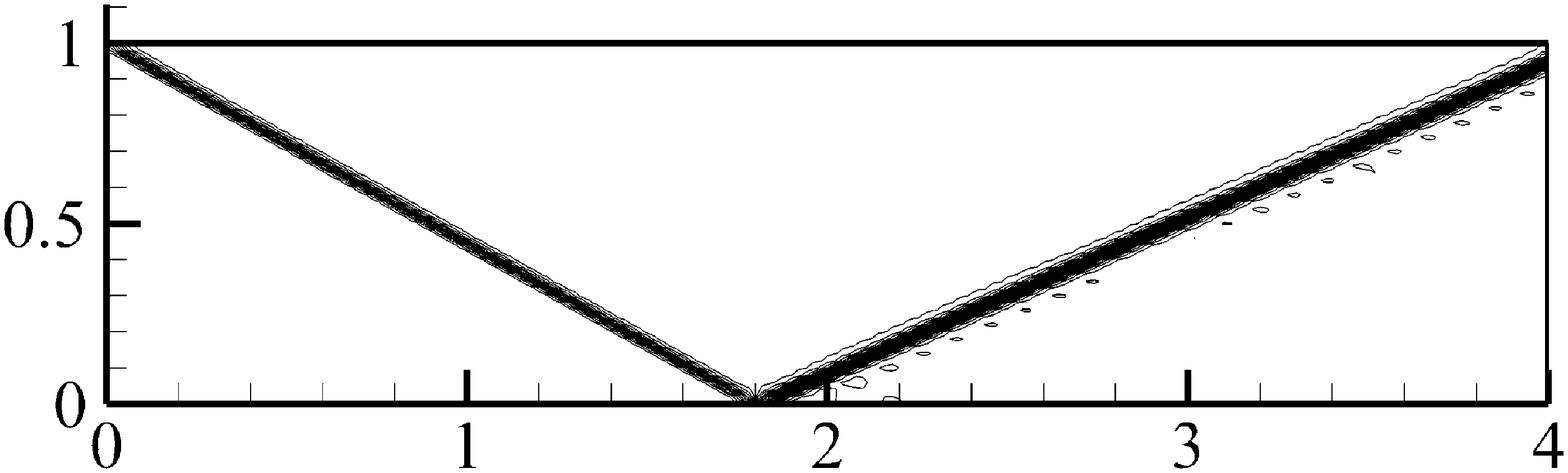}}
  \subfigure[WENO-FV with yCD]{
  \label{FIG:Shock_Ref_yCD}
  \includegraphics[width=10 cm]{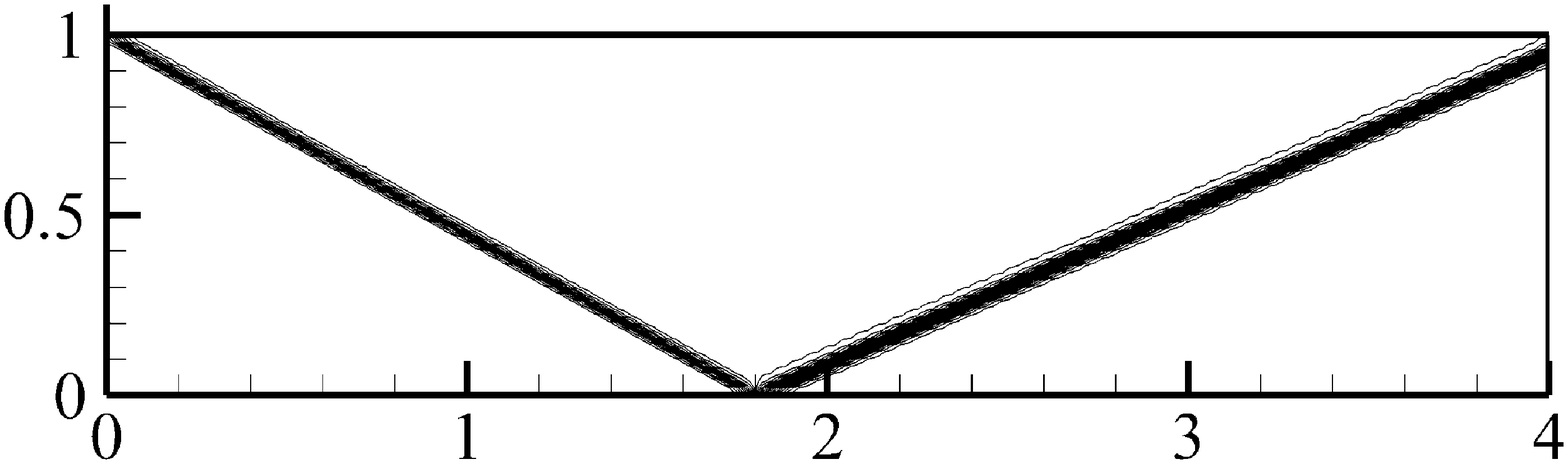}}
  \caption{The density contours of the steady shock reflection on a plate at $t=15$
  calculated by the third-order WENO-FV scheme with different characteristic decomposition techniques.
  The mesh size is 1/50 and the density contours contain 30 equidistant contours from 1.1 to 2.7.}
\label{FIG:Shock_Ref_Density_Contour}
\end{figure}

\begin{figure}
  \centering
  \subfigure[WENO-FV with NCD]{
  \label{FIG:Shock_Ref_Center_NCD}
  \includegraphics[width=10 cm]{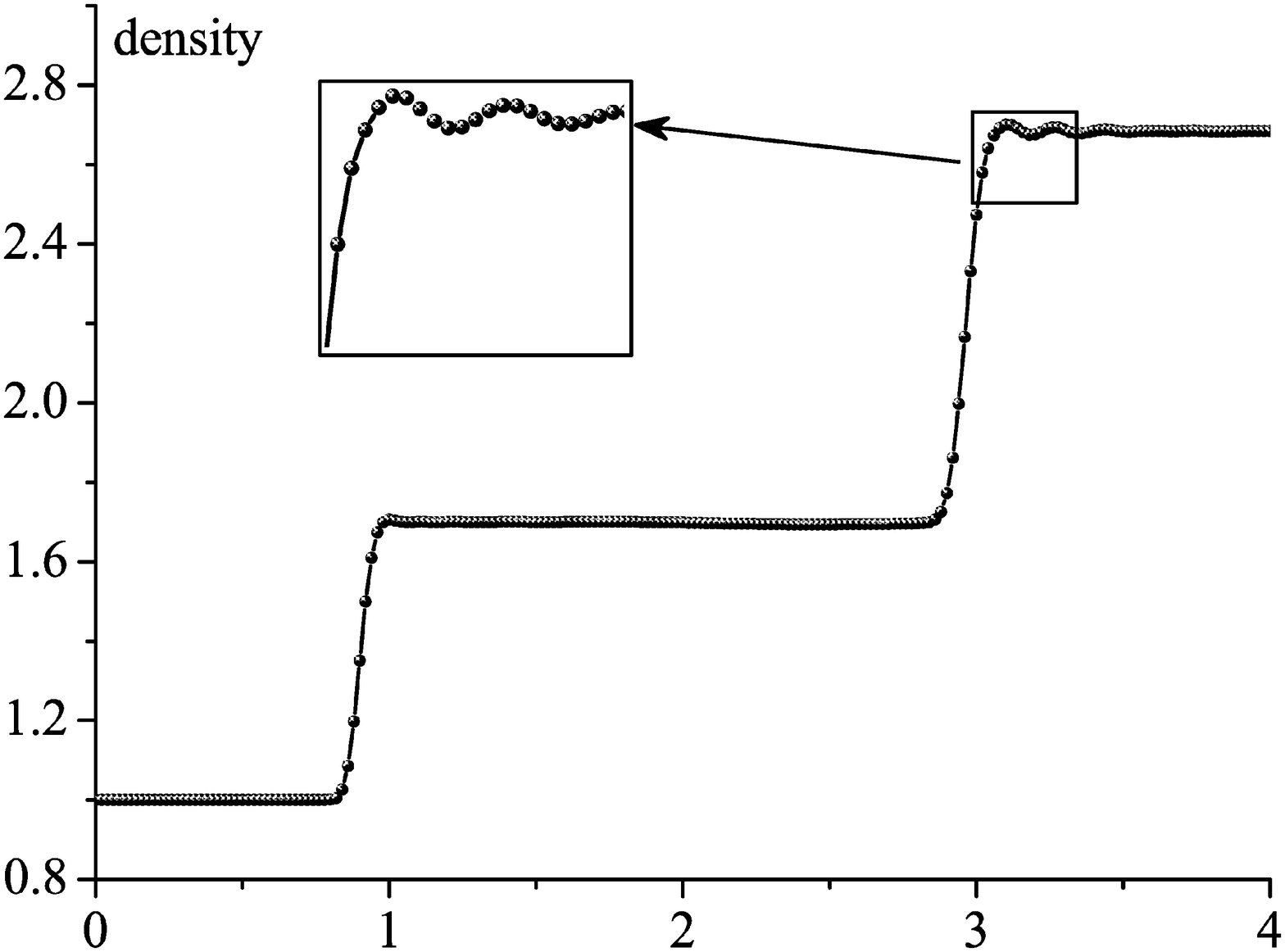}}
  \subfigure[WENO-FV with SCD]{
  \label{FIG:Shock_Ref_Center_SCD}
  \includegraphics[width=10 cm]{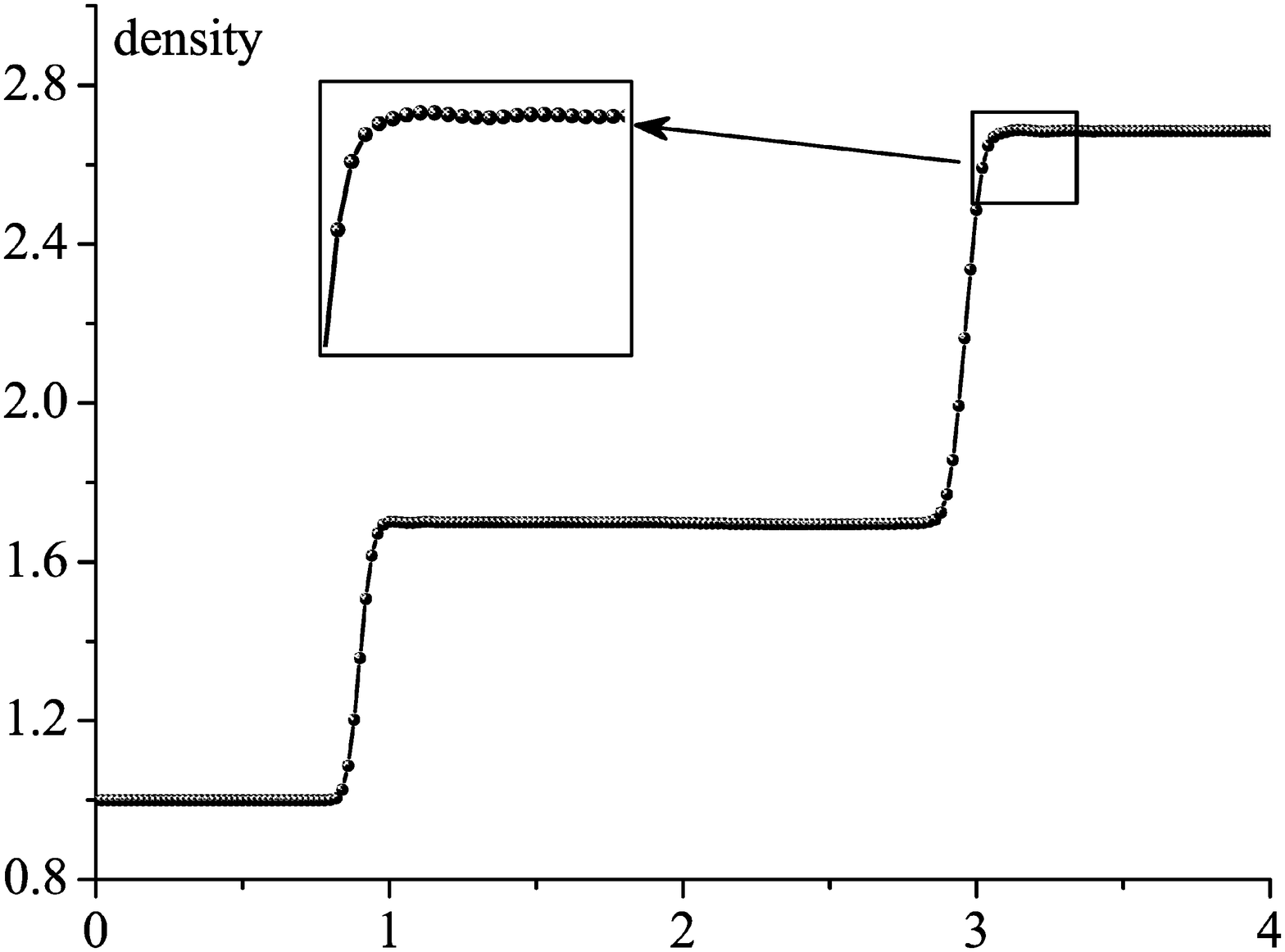}}
  \subfigure[WENO-FV with RCD]{
  \label{FIG:Shock_Ref_Center_RCD}
  \includegraphics[width=10 cm]{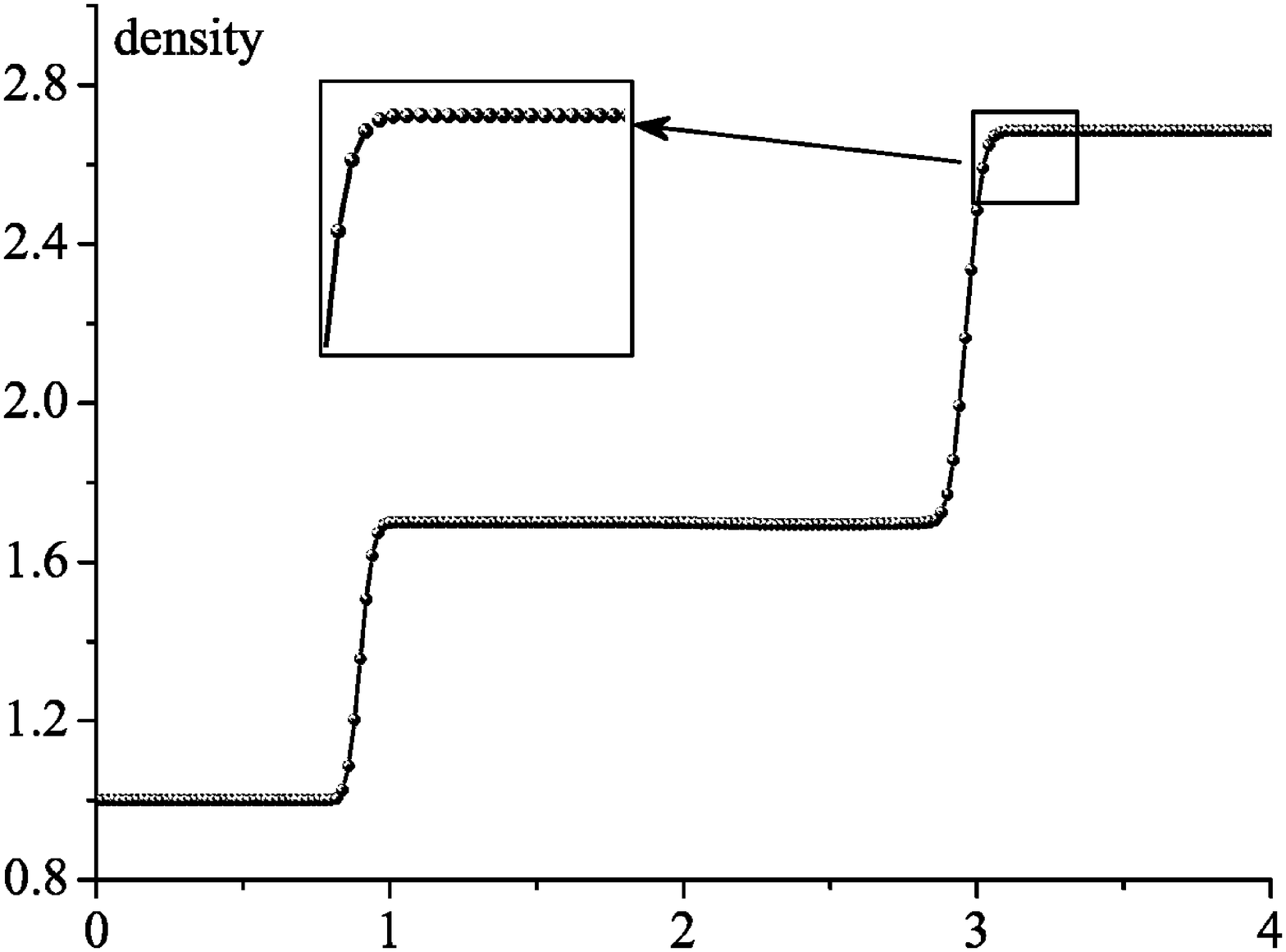}}
  \addtocounter{figure}{-1}      
\end{figure}

\begin{figure} 
  \centering
  \addtocounter{figure}{1}     
  \subfigure[WENO-FV with xCD]{
    \label{FIG:Shock_Ref_Center_xCD}
    \includegraphics[width=10 cm]{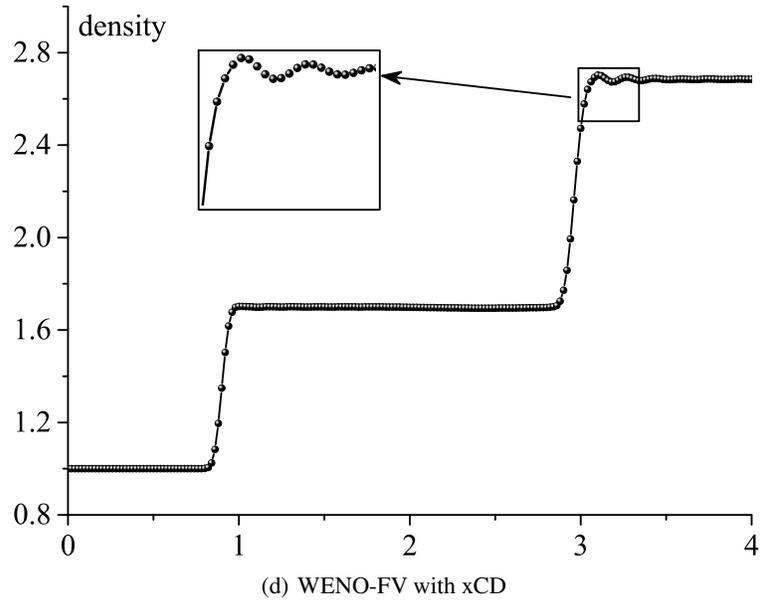}}
    \subfigure[WENO-FV with yCD]{
    \label{FIG:Shock_Ref_Center_yCD}
    \includegraphics[width=10 cm]{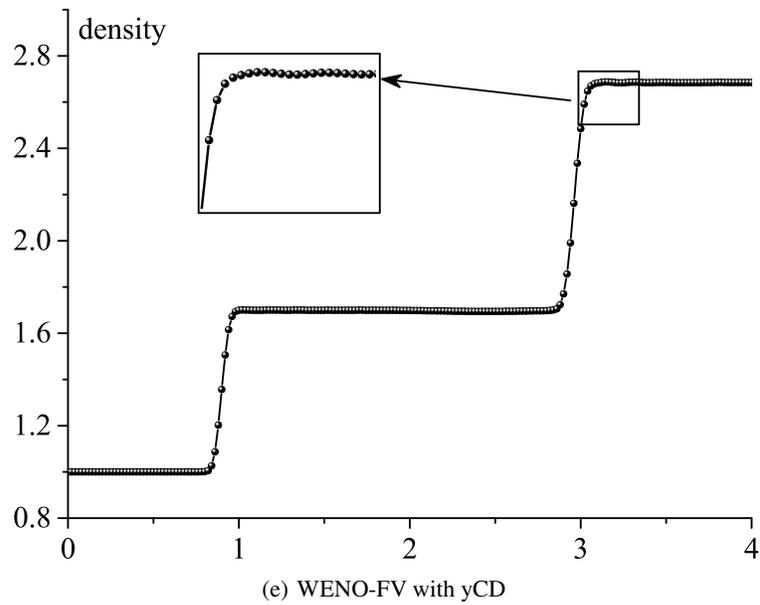}}
  \caption{The density profiles of the steady shock reflection on a plate along the line $y=0.5$ at $t=15$
  calculated by the third-order WENO-FV scheme with different characteristic decomposition techniques.}
\label{FIG:Shock_Ref_Center_Density}
\end{figure}

\begin{figure}
  \centering
  \subfigure[WENO-FV with NCD and SCD]{
  \label{FIG:Shock_Ref_Residual_SCD}
  \includegraphics[width=10 cm]{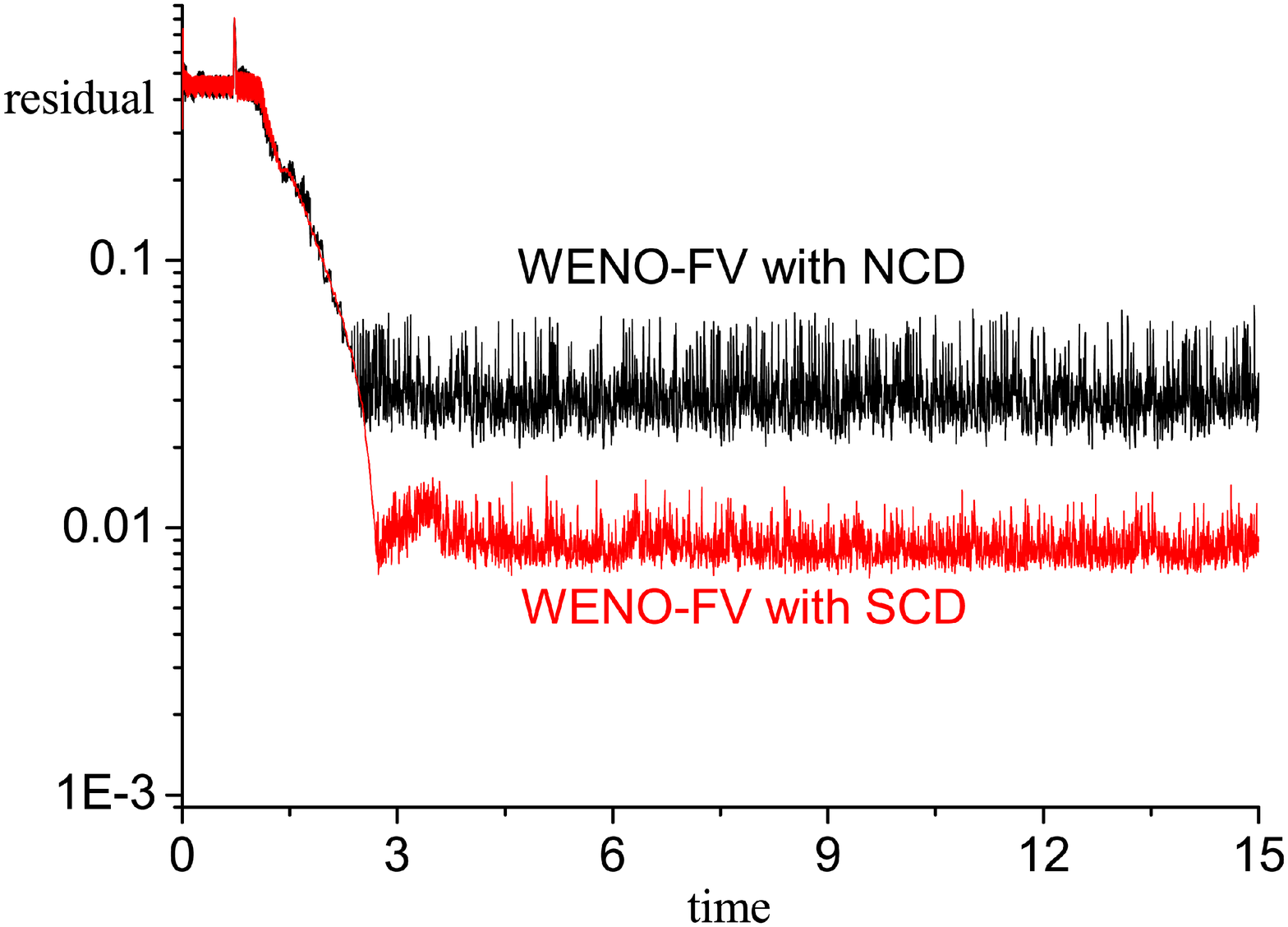}}
  \subfigure[WENO-FV with NCD, xCD, yCD, and RCD]{
  \label{FIG:Shock_Ref_Residual_RCD}
  \includegraphics[width=10 cm]{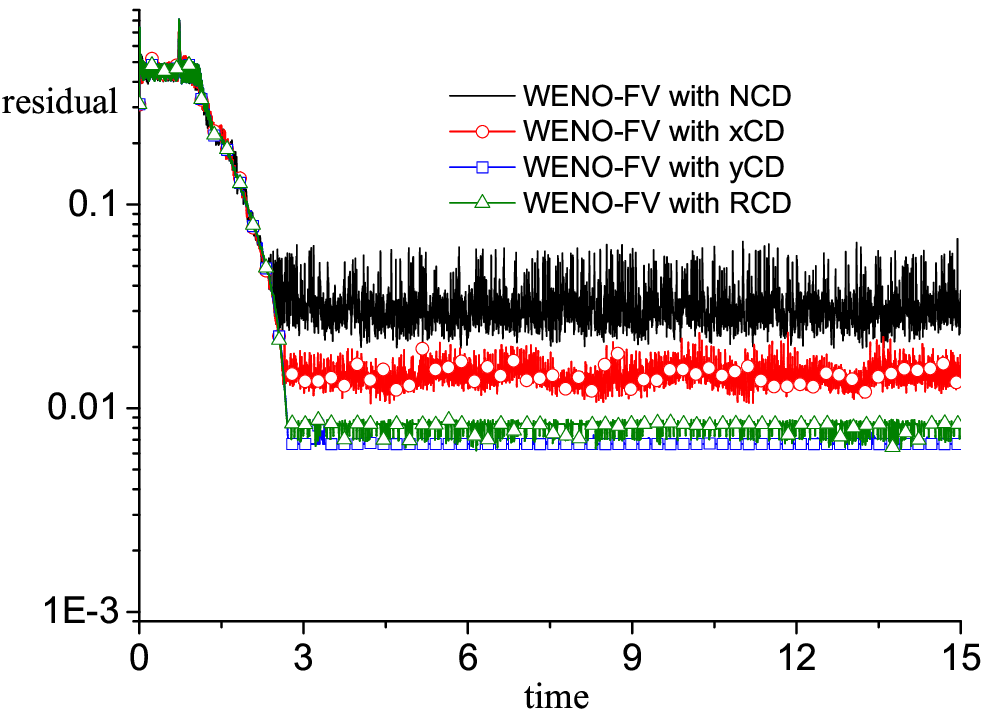}}
  \caption{The maximum residual histories of the steady shock reflection calculated
  by the third-order WENO-FV scheme with different characteristic decomposition techniques.}
\label{FIG:Shock_Ref_Residual}
\end{figure}

We calculate this problem by the third-order WENO-FV scheme with $200\times50$ cells up to $t=15$.
The CPU times of WENO-FV with NCD, SCD, and RCD are 220s, 598s, and 358s, respectively.
Fig. \ref{FIG:Shock_Ref_Density_Contour} shows the density contours,
and Fig. \ref{FIG:Shock_Ref_Center_Density} shows the corresponding density profiles along the line $y=0.5$.
Again, obvious oscillations are observed near the shocks calculated by the WENO-FV scheme with NCD.
It is interesting to observe that RCD controls spurious oscillations better than SCD.
In Fig. \ref{FIG:Shock_Ref_Center_SCD}, slight oscillations can also be observed
in the density profile along the line $y=0.5$ derived by the WENO-FV scheme with SCD.
In contrast, oscillations are invisible in the density profile derived by the WENO-FV scheme with RCD.
In order to confirm the embedded mechanism,
we additionally perform the characteristic decomposition in $x$ direction (xCD) and $y$ direction (yCD).
We observe that yCD achieves better results than xCD,
but slight oscillations near the reflected shock still exist as shown by Fig. \ref{FIG:Shock_Ref_Center_yCD}.
Since $y$ direction is closer to the normal direction of the shock than $x$ direction,
we conclude that the good outcome of RCD comes from 
the fact that the shock is oblique with respect to the grid and RCD works in the direction normal to it.
Furthermore, we compare the maximum residual histories
of the third-order WENO-FV scheme with different characteristic decomposition techniques in Fig. \ref{FIG:Shock_Ref_Residual}.
Here, the maximum residual is calculated by
$R_{max}^n=MAX(R_1^n,R_2^n,R_3^n,R_4^n)$ where $R_k^n=MAX \left|U_k^n-U_k^{n-1}\right|$ $(k=1\text{ to }4)$ for all cells.
As we can see from Fig. \ref{FIG:Shock_Ref_Residual}, yCD has the
smallest residuals among all the techniques, and RCD comes second. 
In addition, yCD and RCD obviously have smaller residual oscillations
than the other three techniques.
We note that, the residuals cannot achieve machine zero regardless of the characteristic decomposition technique.
This is caused by the WENO reconstruction technique.
If we want to reduce the residuals to machine zero, we have to use other technique,
see for example the work of Zhu and Shu \cite{Zhu2017Steady}.

\section{Conclusions}
We propose a rotated characteristic decomposition technique for high-order reconstructions in multi-dimensions.
The third-order WENO-FV scheme with the proposed technique for the Euler equations
can eliminate spurious effectively with about $40\%$ less CPU times than the standard characteristic decomposition technique.
For the steady shock reflection problem, the rotated characteristic decomposition technique achieves a even better performance
than the standard characteristic decomposition technique does.
The proposed technique can apply to other hyperbolic systems and other numerical methods relying on reconstructions.
It is also easy to extend to three-dimensional and unstructured computations.
Moreover, since the rotated characteristic decomposition technique requires only one-time characteristic decomposition 
regardless the number of normal directions of the mesh, 
it can save more CPU times for the three-dimensional and unstructured computations than 
the current Cartesian computations in two dimensions.
If we want to further reduce the computational cost,
we can employ the current approach in addition to the existing ones
that only perform characteristic decomposition near discontinuities
\cite{Ren2003hybridCompatWENO, Puppo2003adaptive, Puppo2011adaptive, Li2010hybridWENO, Peng2019adaptiveWENOZ}.

\bibliographystyle{unsrtnat}
\bibliography{references}  

\begin{thebibliography}{38}
\providecommand{\natexlab}[1]{#1}
\providecommand{\url}[1]{\texttt{#1}}
\expandafter\ifx\csname urlstyle\endcsname\relax
  \providecommand{\doi}[1]{doi: #1}\else
  \providecommand{\doi}{doi: \begingroup \urlstyle{rm}\Url}\fi

\bibitem[Van~Leer(1979)]{vanLeer1979MUSCLV}
Bram Van~Leer.
\newblock Towards the ultimate conservative difference scheme. {V}. {A}
  second-order sequel to {G}odunov's method.
\newblock \emph{Journal of computational Physics}, 32\penalty0 (1):\penalty0
  101--136, 1979.

\bibitem[Godunov(1959)]{Godunov1959}
Sergei~Konstantinovich Godunov.
\newblock A difference method for numerical calculation of discontinuous
  solutions of the equations of hydrodynamics.
\newblock \emph{Matematicheskii Sbornik}, 89\penalty0 (3):\penalty0 271--306,
  1959.

\bibitem[Colella and Woodward(1984)]{Colella1984PPM}
Phillip Colella and Paul~R Woodward.
\newblock The piecewise parabolic method ({PPM}) for gas-dynamical simulations.
\newblock \emph{Journal of computational physics}, 54\penalty0 (1):\penalty0
  174--201, 1984.

\bibitem[Harten et~al.(1987)Harten, Engquist, Osher, and
  Chakravarthy]{Harten1987ENO}
Ami Harten, Bjorn Engquist, Stanley Osher, and Sukumar~R Chakravarthy.
\newblock Uniformly high order accurate essentially non-oscillatory schemes,
  {III}.
\newblock \emph{Journal of Computational Physics}, 71\penalty0 (2):\penalty0
  231--303, 1987.
\newblock \doi{10.1016/0021-9991(87)90031-3}.

\bibitem[Liu et~al.(1994)Liu, Osher, Chan, et~al.]{Liu1994WENO}
Xu-Dong Liu, Stanley Osher, Tony Chan, et~al.
\newblock Weighted essentially non-oscillatory schemes.
\newblock \emph{Journal of computational physics}, 115\penalty0 (1):\penalty0
  200--212, 1994.

\bibitem[Jiang and Shu(1996)]{Jiang_Shu1996WENO}
Guang-Shan Jiang and Chi-Wang Shu.
\newblock Efficient implementation of weighted {ENO} schemes.
\newblock \emph{Journal of Computational Physics}, 126\penalty0 (1):\penalty0
  202--228, 1996.
\newblock \doi{10.1006/jcph.1996.0130}.

\bibitem[Borges et~al.(2008)Borges, Carmona, Costa, and Don]{Borges2008WENO_Z}
Rafael Borges, Monique Carmona, Bruno Costa, and Wai~Sun Don.
\newblock An improved weighted essentially non-oscillatory scheme for
  hyperbolic conservation laws.
\newblock \emph{Journal of Computational Physics}, 227\penalty0 (6):\penalty0
  3191--3211, 2008.

\bibitem[Levy et~al.(1999)Levy, Puppo, and Russo]{Levy1999CWENO}
Doron Levy, Gabriella Puppo, and Giovanni Russo.
\newblock Central {WENO} schemes for hyperbolic systems of conservation laws.
\newblock \emph{ESAIM: Mathematical Modelling and Numerical Analysis},
  33\penalty0 (3):\penalty0 547--571, 1999.

\bibitem[Levy et~al.(2000{\natexlab{a}})Levy, Puppo, and
  Russo]{Levy2000CWENO_ANM}
Doron Levy, Gabriella Puppo, and Giovanni Russo.
\newblock A third order central {WENO} scheme for 2{D} conservation laws.
\newblock \emph{Applied Numerical Mathematics}, 33\penalty0 (1):\penalty0
  415--422, 2000{\natexlab{a}}.

\bibitem[Levy et~al.(2000{\natexlab{b}})Levy, Puppo, and
  Russo]{Levy2000CompactCWENO}
Doron Levy, Gabriella Puppo, and Giovanni Russo.
\newblock Compact central {WENO} schemes for multidimensional conservation
  laws.
\newblock \emph{SIAM Journal on Scientific Computing}, 22\penalty0
  (2):\penalty0 656--672, 2000{\natexlab{b}}.

\bibitem[Cravero et~al.(2019)Cravero, Semplice, and Visconti]{Cravero2019CWENO}
Isabella Cravero, Matteo Semplice, and Giuseppe Visconti.
\newblock Optimal definition of the nonlinear weights in multidimensional
  central {WENOZ} reconstructions.
\newblock \emph{SIAM Journal on Numerical Analysis}, 57\penalty0 (5):\penalty0
  2328--2358, 2019.

\bibitem[Abgrall(1994)]{Abgrall1994ENO_Tri}
R{\'e}mi Abgrall.
\newblock On essentially non-oscillatory schemes on unstructured meshes:
  analysis and implementation.
\newblock \emph{Journal of Computational Physics}, 114\penalty0 (1):\penalty0
  45--58, 1994.

\bibitem[Hu and Shu(1999)]{Hu1999WENO_Tri}
Changqing Hu and Chi-Wang Shu.
\newblock Weighted essentially non-oscillatory schemes on triangular meshes.
\newblock \emph{Journal of Computational Physics}, 150\penalty0 (1):\penalty0
  97--127, 1999.

\bibitem[Dumbser and K{\"a}ser(2007)]{Dumbser2007WENO_FV}
Michael Dumbser and Martin K{\"a}ser.
\newblock Arbitrary high order non-oscillatory finite volume schemes on
  unstructured meshes for linear hyperbolic systems.
\newblock \emph{Journal of Computational Physics}, 221\penalty0 (2):\penalty0
  693--723, 2007.

\bibitem[Reed and Hill(1973)]{Reed1973DG}
Wm~H Reed and TR~Hill.
\newblock Triangular mesh methods for the neutron transport equation.
\newblock \emph{Los Alamos Report LA-UR-73-479}, 1973.

\bibitem[Cockburn and Shu(1989)]{Cockburn1989RKDGII}
Bernardo Cockburn and Chi-Wang Shu.
\newblock {TVB} {R}unge-{K}utta local projection discontinuous {G}alerkin
  finite element method for conservation laws. {II}. {G}eneral framework.
\newblock \emph{Mathematics of computation}, 52\penalty0 (186):\penalty0
  411--435, 1989.

\bibitem[Cockburn et~al.(1989)Cockburn, Lin, and Shu]{Cockburn1989RKDGIII}
Bernardo Cockburn, San-Yih Lin, and Chi-Wang Shu.
\newblock {TVB} {R}unge-{K}utta local projection discontinuous {G}alerkin
  finite element method for conservation laws. {III}: {O}ne-dimensional
  systems.
\newblock \emph{Journal of computational physics}, 84\penalty0 (1):\penalty0
  90--113, 1989.

\bibitem[Cockburn et~al.(1990)Cockburn, Hou, and Shu]{Cockburn1990RKDGIV}
Bernardo Cockburn, Suchung Hou, and Chi-Wang Shu.
\newblock The {R}unge-{K}utta local projection discontinuous {G}alerkin finite
  element method for conservation laws. {IV}. {T}he multidimensional case.
\newblock \emph{Mathematics of Computation}, 54\penalty0 (190):\penalty0
  545--581, 1990.

\bibitem[Cockburn and Shu(1998)]{Cockburn1998RKDGV}
Bernardo Cockburn and Chi-Wang Shu.
\newblock The {R}unge-{K}utta local projection discontinuous {G}alerkin finite
  element method for conservation laws. {V}: {M}ultidimensional systems.
\newblock \emph{Journal of computational physics}, 141\penalty0 (2):\penalty0
  199--224, 1998.

\bibitem[Qiu and Shu(2004)]{Qiu2004HermiteWENO}
Jianxian Qiu and Chi-Wang Shu.
\newblock Hermite {WENO} schemes and their application as limiters for
  {R}unge--{K}utta discontinuous {G}alerkin method: {O}ne-dimensional case.
\newblock \emph{Journal of Computational Physics}, 193\penalty0 (1):\penalty0
  115--135, 2004.

\bibitem[Qiu and Shu(2005)]{Qiu2005HermiteWENO2D}
Jianxian Qiu and Chi-Wang Shu.
\newblock Hermite {WENO} schemes and their application as limiters for
  {R}unge--{K}utta discontinuous {G}alerkin method {II}: {T}wo-dimensional
  case.
\newblock \emph{Computers \& Fluids}, 34\penalty0 (6):\penalty0 642--663, 2005.

\bibitem[Zhong and Shu(2013)]{Zhong2013DG_Limiter}
Xinghui Zhong and Chi-Wang Shu.
\newblock A simple weighted essentially nonoscillatory limiter for
  {R}unge--{K}utta discontinuous {G}alerkin methods.
\newblock \emph{Journal of Computational Physics}, 232\penalty0 (1):\penalty0
  397--415, 2013.

\bibitem[Zhu et~al.(2016)Zhu, Zhong, Shu, and Qiu]{Zhu2016DG_Limiter}
Jun Zhu, Xinghui Zhong, Chi-Wang Shu, and Jianxian Qiu.
\newblock Runge-{K}utta discontinuous {G}alerkin method with a simple and
  compact {H}ermite {WENO} limiter.
\newblock \emph{Communications in Computational Physics}, 19\penalty0
  (4):\penalty0 944--969, 2016.

\bibitem[Dumbser et~al.(2008)Dumbser, Balsara, Toro, and Munz]{Dumbser2008PNPM}
Michael Dumbser, Dinshaw~S Balsara, Eleuterio~F Toro, and Claus-Dieter Munz.
\newblock A unified framework for the construction of one-step finite volume
  and discontinuous {G}alerkin schemes on unstructured meshes.
\newblock \emph{Journal of Computational Physics}, 227\penalty0 (18):\penalty0
  8209--8253, 2008.
\newblock \doi{10.1016/j.jcp.2008.05.025}.

\bibitem[Dumbser and Zanotti(2009)]{Dumbser2009PNPM}
Michael Dumbser and Olindo Zanotti.
\newblock Very high order {$P_NP_M$} schemes on unstructured meshes for the
  resistive relativistic {MHD} equations.
\newblock \emph{Journal of Computational Physics}, 228\penalty0 (18):\penalty0
  6991--7006, 2009.

\bibitem[Dumbser(2010)]{Dumbser2010PNPM}
Michael Dumbser.
\newblock Arbitrary high order {$P_NP_M$} schemes on unstructured meshes for
  the compressible {N}avier--{S}tokes equations.
\newblock \emph{Computers \& Fluids}, 39\penalty0 (1):\penalty0 60--76, 2010.

\bibitem[Qiu and Shu(2002)]{Qiu2002Decomposition}
Jianxian Qiu and Chi-Wang Shu.
\newblock On the construction, comparison, and local characteristic
  decomposition for high-order central weno schemes.
\newblock \emph{Journal of Computational Physics}, 183\penalty0 (1):\penalty0
  187--209, 2002.

\bibitem[Shu and Osher(1988)]{Shu1988EfficientENO}
Chi-Wang Shu and Stanley Osher.
\newblock Efficient implementation of essentially non-oscillatory
  shock-capturing schemes.
\newblock \emph{Journal of Computational Physics}, 77:\penalty0 439--471, 1988.

\bibitem[Ren et~al.(2003)Ren, Zhang, et~al.]{Ren2003hybridCompatWENO}
Yu-Xin Ren, Hanxin Zhang, et~al.
\newblock A characteristic-wise hybrid compact-weno scheme for solving
  hyperbolic conservation laws.
\newblock \emph{Journal of Computational Physics}, 192\penalty0 (2):\penalty0
  365--386, 2003.

\bibitem[Puppo(2003)]{Puppo2003adaptive}
Gabriella Puppo.
\newblock Adaptive application of characteristic projection for central
  schemes.
\newblock In \emph{Hyperbolic problems: theory, numerics, applications}, pages
  819--829. Springer, 2003.

\bibitem[Puppo and Semplice(2011)]{Puppo2011adaptive}
Gabriella Puppo and Matteo Semplice.
\newblock Numerical entropy and adaptivity for finite volume schemes.
\newblock \emph{Communications in Computational Physics}, 10\penalty0
  (5):\penalty0 1132--1160, 2011.

\bibitem[Li and Qiu(2010)]{Li2010hybridWENO}
Gang Li and Jianxian Qiu.
\newblock Hybrid weighted essentially non-oscillatory schemes with different
  indicators.
\newblock \emph{Journal of Computational Physics}, 229\penalty0 (21):\penalty0
  8105--8129, 2010.

\bibitem[Peng et~al.(2019)Peng, Zhai, Ni, Yong, and
  Shen]{Peng2019adaptiveWENOZ}
Jun Peng, Chuanlei Zhai, Guoxi Ni, Heng Yong, and Yiqing Shen.
\newblock An adaptive characteristic-wise reconstruction weno-z scheme for gas
  dynamic euler equations.
\newblock \emph{Computers \& Fluids}, 179:\penalty0 34--51, 2019.

\bibitem[Balsara et~al.(2009)Balsara, Rumpf, Dumbser, and
  Munz]{Balsara2009ADER-WENO}
Dinshaw~S Balsara, Tobias Rumpf, Michael Dumbser, and Claus-Dieter Munz.
\newblock Efficient, high accuracy {ADER-WENO} schemes for hydrodynamics and
  divergence-free magnetohydrodynamics.
\newblock \emph{Journal of Computational Physics}, 228\penalty0 (7):\penalty0
  2480--2516, 2009.

\bibitem[Toro(2013)]{Toro2013CFDBook}
Eleuterio~F Toro.
\newblock \emph{Riemann solvers and numerical methods for fluid dynamics: a
  practical introduction}.
\newblock Springer Science \& Business Media, 2013.

\bibitem[Woodward and Colella(1984)]{Woodward1984JCP}
Paul Woodward and Phillip Colella.
\newblock The numerical simulation of two-dimensional fluid flow with strong
  shocks.
\newblock \emph{Journal of computational physics}, 54\penalty0 (1):\penalty0
  115--173, 1984.

\bibitem[Yee et~al.(1985)Yee, Warming, and Harten]{Yee1985JCP}
HC~Yee, RF~Warming, and A~Harten.
\newblock Implicit total variation diminishing {(TVD)} schemes for steady-state
  calculations.
\newblock \emph{Journal of Computational Physics}, 57\penalty0 (3):\penalty0
  327--360, 1985.

\bibitem[Zhu and Shu(2017)]{Zhu2017Steady}
Jun Zhu and Chi-Wang Shu.
\newblock Numerical study on the convergence to steady state solutions of a new
  class of high order {WENO} schemes.
\newblock \emph{Journal of Computational Physics}, 349:\penalty0 80--96, 2017.

\end{thebibliography}






\end{document}